\documentclass[10pt]{article}
\usepackage{graphicx} % Required for inserting images
\usepackage{amsmath,amssymb,enumerate,xcolor}
\usepackage[hidelinks]{hyperref}

\title{On the Jordan--Chevalley decomposition problem for operator fields in small dimensions and Tempesta--Tondo  conjecture}
\author{Alexey V.\ Bolsinov\footnote{School of Mathematics,
 Loughborough University,
 LE11 3TU, UK  and Institute of Mathematics and Mathematical Modeling, Almaty, Kazakhstan\ \ \quad {\tt A.Bolsinov@lboro.ac.uk}},    
Andrey Yu.\  Konyaev\footnote{Faculty of Mechanics and Mathematics and Center for Fundamental and Applied Mathematics, Moscow State University, 119992, Moscow, Russia  and Institute of Mathematics and Mathematical Modeling, Almaty, Kazakhstan
 \ \ \quad{\tt  maodzund@yandex.ru}}  \, and 
    Vladimir S.\ Matveev\footnote{
Institut f\"ur Mathematik, Friedrich Schiller Universit\"at Jena,
07737 Jena,  Germany  \ \ \quad {\tt  vladimir.matveev@uni-jena.de}}}

\date{March  2025}

\newtheorem{theorem}{Theorem}
\newtheorem{remark}[theorem]{Remark}

\newtheorem{example}{Example}[section]

\newcommand{\weg}[1]{}

\usepackage[backend=biber]{biblatex}  %Imports biblatex package
\begin{document}

\maketitle
\begin{abstract} 
We explore the Jordan--Chevalley decomposition problem for an operator field in small dimensions. In dimensions three and four,  we find tensorial conditions for an operator field $L$, similar to a nilpotent Jordan block, to possess local coordinates in which $L$ takes a strictly upper triangular form. We prove the Tempesta--Tondo conjecture for higher order brackets of Fr\"olicher-Nijenhuis type.
\end{abstract}

\section{Introduction}

Let $L$ be an operator field, i.e, a tensor field of type $(1, 1)$. The {\it Nijenhuis torsion} of $L$ is the tensor of type $(1, 2)$ given by 
\begin{equation}\label{nij}
   \mathcal T_L (\xi, \eta) =  L^2[\xi, \eta] + [L\xi, L\eta] - L[L\xi, \eta] - L[\xi, L\eta].  
\end{equation}
The next recursion formula defines {\it Haantjes torsion} of level $m$ (see \cite{Kosmann} and \cite{TT22}):
\begin{equation}\label{haam}
\begin{aligned}
\mathcal T^{(1)}_L(\xi, \eta) = & \mathcal T_L(\xi, \eta), \\
\mathcal T^{(m)}_L(\xi, \eta) = & L^2 \mathcal T^{(m - 1)}_L(\xi, \eta) + \mathcal T^{(m - 1)}_L(L\xi, L\eta) - \\
& - L\mathcal T^{(m - 1)}_L(L\xi, \eta) - L\mathcal T^{(m - 1)}_L(\xi, L\eta), \quad m = 2, 3, \dots.   
\end{aligned}    
\end{equation}
For $m=2$, this formula coincides with the classical definition of Haantjes torsion $\mathcal H_L = \mathcal T^{(2)}_L$  introduced in \cite{nijenhuis1951,Nijenhuis1955, haantjes}.

Both Nijenhuis and Haantjes torsions are widely used in many areas of mathematics. To understand the reason for that, recall the famous Haantjes\footnote{Under the additional  assumption that the spectrum of $L$ is simple, the result is due to Nijenhuis \cite{nijenhuis1951}} criterion \cite{haantjes}: 

{\it  For a semi-simple operator field $L$ with real spectrum,   the  Haantjes torsion $\mathcal H_L$ vanishes if and only if in a neighborhood of almost every point  there exists  a local coordinate system  in which $L$ is given by a  diagonal matrix.}

This, differential-geometric by its nature, result is widely used in mathematical physics, in particular, in the theory of evolutionary PDEs of hydrodynamic type \cite{Bo2006, FM2007, Sharipov96} and in the theory of separation of variables for finite-dimensional integrable systems \cite{BKM22,  KKM,  NTT,TT21, TT22b, TT16}. It provides an invariant and calculable condition for the existence of a `good' coordinate system for generic operator fields $L$. In many situations of interest, vanishing of the Haantjes torsion follows from some additional assumptions, and then the Haantjes criterion provides a natural ansatz for $L$,  suitable for further computations. 

The higher Haantjes torsions appeared quite recently and were  independently introduced in \cite[\S 4.4]{Kosmann} and \cite[\S 4.4]{TT22} in the context  of integrable systems. In \cite[Corollary 27]{TT22} it was shown that in dimension $n$,  the Haantjes torsion $\mathcal T^{(n - 1)}_L$ vanishes for any operator field $L$ given by a strictly upper triangular matrix. In other words, the condition $\mathcal T^{(n - 1)}_L = 0$ is necessary for a nilpotent operator field $L$ to be brought to a triangular form.

So we pose a natural question, partially motivated by the discussion in \cite[\S\S 4.2, 4.3]{TT22}.  Let $L$ be similar to a $n\times n$ nilpotent Jordan block at every point. Does the vanishing of $\mathcal T^{(n - 1)}_L$ provide a sufficient condition for the existence of coordinates, in which $L$ is strictly upper triangular? The following example shows that for $n \geq 4$ the answer is negative.

\begin{example} \label{ex:1}
\rm{ Let $L$ be an operator field, which is similar to a nilpotent Jordan block at each point.  Assume that $L$ can be brought to a (strictly) upper triangular form.  Then one can easily see that all the distributions in the flag
$$
\{ 0\} \subset \operatorname{Image}  L ^{n-1}\subset \dots \subset \operatorname{Image}  L 
$$
are integrable. The converse is also true.
Consider now the operator field  given by the matrix 
\begin{equation}\label{eq:ex1}
L =
\begin{pmatrix}
0 & 1 & 0 & 0 \\
0 & 0 & 1 & 0 \\
0 & 0 & -x_2 & 1 \\
0 & 0 & -x_2^2 & x_2
\end{pmatrix}.
\end{equation}
This operator is nilpotent and similar to a Jordan block at each point. The image of  $L$ is spanned by three vector fields
$$
\xi_1 =  \partial_{x_1},   \quad \xi_2 = \partial_{x_2}, \quad \xi_3= \partial_{x_3} + x_2\partial_{x_4}.
$$
Notice that $[\xi_2, \xi_3] = \partial_{x_4} \notin \operatorname{span}(\xi_1, \xi_2,\xi_3)$. Hence,  the  distribution $\operatorname{Image} L$
is not integrable and, therefore, $L$ cannot be put to a strictly upper triangular form by a coordinate change. It is straightforward to check, however, that  $\mathcal H^{(3)}_L\equiv 0$. 
}
\end{example}

Let us briefly discuss dimensions 2 and 3. In dimension 2, the Haantjes torsion of every operator $L$ vanishes and, according to the Haantjes criterion, every operator $L$ having two different real eigenvalues can be locally diagonalized. 

Now assume that $L$ has only one real eigenvalue $f$ and $\operatorname{rank}(L- f \,\textrm{\bf 1}) = 1\ne 0$ so that $L$ is not diagonalisable at any point.  Then locally there exists a (unique up to proportionality) smooth vector field $\xi$ such that $L\xi = f\,\xi$.  If we choose local coordinates $x, y$ in such a way that  $\xi = \partial_{x}$, then $L$ automatically takes the triangular form 
$$
 L = \begin{pmatrix}
     f(x, y) & g(x, y)  \\
     0 & f(x, y)
\end{pmatrix}.
$$

Thus, in dimension 2, local reduction to a good normal form requires no additional conditions.  
Recall, however, that in the context of diagonalisability and/or triangularisability problem for operator fields $L$, one still needs to assume that the multiplicities of eigenvalues $f_i$ of $L$, as well as the ranks of $(L- f_i \,\textrm{\bf 1})^k$ remain locally constant.   Indeed, the operator 
$A=\begin{pmatrix} 2x & y \\ y & 0 \end{pmatrix}$ is $\mathbb R$-diagonalisable pointwise and the Nijenhuis torsion $\mathcal T_A$ vanishes. However, in a neighbourhood of the point $\mathsf p = (0,0)$, it can be reduced to neither diagonal nor triangular form.  The reason is that the eigenvalues of $L$ collide at this point. Another example is  $B = \begin{pmatrix}  xy & -y^2 \\ x^2 & -xy  \end{pmatrix}$.  This operator is nilpotent and $\mathcal T_B=0$, but $\operatorname{rank}L$ drops at  $\mathsf p=(0,0)$ and, as a result, $B$ is not reducible to a triangular form in any neighbourhood of $\mathsf p$.

In the case of dimension 3, the following theorem holds.

\begin{theorem} \label{thm:dim3a} 
Let $L$ be an operator field in dimension three such that at every point $L$ has only one eigenvalue, and this eigenvalue has geometric multiplicity one. Then the following are equivalent:
\begin{itemize}
\item in a neighbourhood of each point, there exists a coordinate system such that  $L$ is upper triangular, 
\item the Haantjes torsion  $\mathcal H_L$  of $L$  vanishes.
\end{itemize}  
\end{theorem}

\begin{example}
   Consider the operator $L$  field given,  in a local coordinate system  $x_1, x_2, x_3$,  by the matrix 
$$\begin{pmatrix}
44 x_1 ^2  - 16 x_1  x_2  + 43 x_2  + 45 x_3 &  66 x_1 ^2  - 20 x_1  x_2  + 66 x_2  + 66 x_3 &       
       55 x_1 ^2 - 24 x_1  x_2  + 55 x_2  + 55 x_3  \\ 
       -16 x_1 ^2 + 8 x_1  x_2  - 16 x_2  - 16 x_3 & 
       -24 x_1 ^2  + 10 x_1  x_2  - 25 x_2  - 23 x_3 & 
       -20 x_1 ^2  + 12 x_1  x_2  - 20 x_2  - 20 x_3 \\  
       -16 x_1 ^2  + 4 x_1  x_2  - 16 x_2  - 16 x_3 &                                                     
       -24 x_1 ^2  + 5 x_1  x_2  - 24 x_2  - 24 x_3 &  
       -20 x_1 ^2  + 6 x_1  x_2  - 21 x_2  - 19 x_3 \end{pmatrix}
. $$ 
The operator is similar to the $3\times 3$-Jordan block with the eigenvalue $x_3 -x_2$. 
The  Haantjes torsion of this operator  is zero, so  there exists a coordinate system such that the matrix of $L$ is upper triangular. 
\end{example}

\begin{example} Consider the operator field  $L$ given,  in a local coordinate system  $x_1, x_2, x_3$,  by the upper diagonal matrix 
$$
\begin{pmatrix} x_1 & x_2 & 0 \\ 0 & x_2 & x_2 \\ 0 & 0 & x_3      
\end{pmatrix}.$$
Its Haantjes torsion is not zero. The example  shows that the assumption in Theorem \ref{thm:dim3a}
that the operator is similar to a single  Jordan block  is essential. 
\end{example}

Now let us focus our attention on dimension 4 and ask the following natural question: is there a tensor field $T$, constructed from an operator  $L$, which vanishes if and only if the `upper triangular' coordinate system exists? The answer is positive and is given in the following Theorem.

\begin{theorem}\label{thm:2}
Let $L$ be an operator field  in dimension four  such that at every point it has only one eigenvalue, and this eigenvalue has geometric multiplicity $1$.
Consider its Haantjes torsion $\mathcal T^{(2)}_L$, which we denote by $\mathcal{ H}^i_{jk}$, and the tensor field 
\begin{equation}
\label{eq:tensorT}
T^i_{jk}= \widehat L^i_s \mathcal{H}^s_{rk} \widehat L^r_j- \widehat L^i_s \mathcal{H}^s_{jr} \widehat L^r_k + \mathcal{H}^i_{sk}\widehat L^s_r \widehat L^r_j
\end{equation}
with $\widehat L:= L- \tfrac{1}{4} \operatorname{trace}L \cdot \textrm{\bf 1}$. Then in a neighbourhood of each point, there exists a coordinate system such that $L$ is upper triangular if and only if $T$ vanishes.
\end{theorem}

The components of the tensor field $T$ from Theorem \ref{thm:2} are polynomials of order five in components $L^i_j$ and linear in derivatives $\frac{\partial L^i_j}{\partial x_p}$. This is an example of the so-called natural differential tensor operation (see \cite{kolar}). In section \ref{how} we provide an algorithm that allows one to search for similar tensors in higher dimensions.

\begin{example} Consider the  operator field  $L$ given, in a local coordinate system, by the upper triangular  matrix
$$\begin{pmatrix}0 & a & 0 & 0\\ 0&0&b &0\\ 0&0&0&c \\0 &0 &0&0\end{pmatrix},$$
where $a,b,c$ are functions of the local coordinates  $x_1,x_2, x_3, x_4.$ 
Then, the entries of the  Haantjes torsion  which can be different from $0$ are:
\begin{eqnarray*}
    \mathcal{H}^1_{42}=-\mathcal{H}^1_{24}&=&2 a^{2} \left(b \left(\frac{\partial c}{\partial  x_{1}}\right)-\left(\frac{\partial  b}{\partial  x_{1}}\right) c\right) \\ 
\mathcal{H}^1_{34} =   -\mathcal{H}^1_{43} &=&  b \left(\left(\frac{\partial  a}{\partial  x_{2}}\right) b c+\left(\frac{\partial  b}{\partial  x_{2}}\right) a c-2 a b \left(\frac{\partial c}{\partial x_{2}}\right)\right)\\
\mathcal{H}^2_{34}=-\mathcal{H}^2_{43} &= &a b \left(b \left(\frac{\partial c}{\partial x_{1}}\right)-\left(\frac{\partial  b}{\partial  x_{1}}\right) c \right). 
\end{eqnarray*}
We see that,  for generic functions   $a,b,c,$  the Haantjes torsion 
is not zero, though  the  operator field  is upper triangular. Clearly, for generic $a,b,c$,  the operator is    similar  to the $4\times 4$ Jordan block with zero eigenvalue,   so the tensor $T$ given by   \eqref{eq:tensorT} vanishes, which can also   be verified by straightforward computation. 
\end{example}

\begin{example}
Consider the  operator field $L$ given, in a local coordinate system, by the following   matrix:
$$\begin{pmatrix}
-x_1 - x_3 - x_2 & -2 x_2 - x_3& -x_1 - 3 x_2 - 2 x_3&   -2 x_2- x_3\\ 
3 x_3 + x_2 + 2 x_1&  
  3 x_3 + 2 x_2 + x_1&  5 x_3 + 4 x_2 + 3 x_1& 
  3 x_3 + 2 x_2 + x_1\\ 
  0&  -x_1 + x_2& -x_1 + x_2& -x_1 + x_2\\ -2 x_3 - x_1& 
  -2 x_3 - 2 x_2 + x_1& -3 x_3 - 3 x_2 &
  -2 x_3 - 2 x_2 + x_1\end{pmatrix}. $$
It is easy to check that,  for generic $x_1, x_2, x_3, x_4$,  the matrix  is similar   to the $4\times 4$ Jordan block with zero eigenvalue.
By direct calculations one sees that the Haantjes torsion is not zero.
However, the tensor $T$ given by  \eqref{eq:tensorT}
vanishes and therefore there exists a coordinate system  in which $L$ is upper triangular.
\end{example}

Note that all the above results can be considered as  special cases of the Jordan--Chevalley decomposition problem for operator fields formulated in \cite[\S  4.2]{TT22}: \noindent{\it Determine under
which conditions there exist coordinate charts  such that an operator field  $L$  can be decomposed into the sum of  two operators,   $L = D + N$, where $D$ is  a diagonal operator and $N$ is a strictly upper triangular operator, commuting with $D$.}

We study the operators in dimensions  three and four in the following  case: the diagonal part $D = f \cdot \textrm{\bf 1}$ for some function $f$ and $N$ is similar to the nilpotent Jordan block of maximal size. In these specific cases,  Theorems \ref{thm:dim3a}  and \ref{thm:2} provide the aforementioned conditions in terms of tensor fields constructed from $L$.

For a pair of operator fields $K, L$, consider the expression
\begin{equation*}
        \begin{aligned}{}
            [[K, L ]](\xi, \eta) & = [K\xi, L\eta] + [L\xi, K\eta] - L[K\xi, \eta] - L[\xi, K\eta] - \\
            & - K[L\xi, \eta] - K[\xi, L\eta] + LK[\xi, \eta] + KL[\xi, \eta].
        \end{aligned}
\end{equation*}
The r.h.s. defines a tensor field of type $(1, 2)$ called {\it Frolicher-Nijenhuis bracket} of operator fields $K$ and $L$. Obviously, $\mathcal T_L = \frac{1}{2} [[L, L]]$. The higher order brackets are given by the recursion formulas
\begin{equation}\label{haam2}
\begin{aligned}
\mathcal H^{(1)}_{K, L}(\xi, \eta) = & [[K, L]](\xi, \eta), \\
\mathcal H^{(m)}_{K, L}(\xi, \eta) = & KL \mathcal H^{(m - 1)}_{K, L}(\xi, \eta) + \mathcal H^{(m - 1)}_{K, L}(K\xi, L\eta) - \\
& - L \mathcal H^{(m - 1)}_{K, L}(K\xi, \eta) - K\mathcal H^{(m - 1)}_{K, L}(\xi, L\eta) + \\
& + LK \mathcal H^{(m - 1)}_{K, L}(\xi, \eta) + \mathcal H^{(m - 1)}_{K, L}(L\xi, K\eta) - \\
& - K \mathcal H^{(m - 1)}_{K, L}(L\xi, \eta) - L\mathcal H^{(m - 1)}_{K, L}(\xi, K\eta), \quad m = 2, 3, \dots \,.   
\end{aligned}    
\end{equation}
These brackets were introduced in \cite{TT22, TT21} as an important tool in the study of Haantjes algebras and their properties related to  integrable systems. As a final result of our work, we prove the conjecture stated in \cite[Conjecture 29]{TT22}: 

\begin{theorem}[Tempesta--Tondo Conjecture] \   \label{thm:1}
{\it  Let $L, M$ be  two commuting, in the algebraic sense,  operators such that in a local coordinate system $x_1,\dots,x_n$ they are both given by strictly upper triangular matrices. 
Then the generalized
Haantjes bracket $\mathcal H^{(n - 1)}_{K, L}$ of level $(n -1)$ vanishes. }
\end{theorem} 

\subsubsection*{Acknowledgements. } A.\,B. and A.\,K. were supported by the Ministry of Science and Higher Education of the Republic of Kazakhstan (grant No. AP23483476). V.\,M.\ was supported by the DFG (projects 455806247 and 529233771), and  by the ARC  (Discovery Programme DP210100951). Part of the results were obtained during a research visit of A.\,B. and A.\,K.  to the Friedrich Schiller University Jena and V.\,M. to Loughborough University supported by the DFG and by the Friedrich Schiller University Jena.
The paper was finalized during the visit of A.\,B. and V.\,M. at the Sydney Mathematical Research Institute (SMRI). We are grateful to SMRI for their hospitality and financial support. 

%%%%%%%%
%%%%%%%%
%%%%%%%%
%%%%%%%%
%%%%%%%%
%%%%%%%%
%%%%%%%%
%%%%%%%%
%%%%%%%%
%%%%%%%%
%%%%%%%%
%%%%%%%%

\section{Proof of Theorems \ref{thm:dim3a} and \ref{thm:2}}

We will  reduce   the  Jordan--Chevalley decomposition problem,  for operators similar to a Jordan block,  to  a linear  algebra problem  which can be handled by computer algebra software and solved by hand in small dimensions 3 and 4. 

We denote the  standard Jordan block of dimension  $n$  with eigenvalue $\lambda$ by $J_n( \lambda), $ for example,
\begin{equation}\label{eq:canon}
J_4(\lambda) =
\begin{bmatrix}
\lambda & 1 & 0 & 0 \\
0 & \lambda & 1 & 0 \\
0 & 0 & \lambda & 1 \\
0 & 0 & 0 & \lambda
\end{bmatrix}.
\end{equation}
We assume that our operator $L$ is similar, at every point $x$, to  $J_n\bigl(\lambda(x)\bigr)$ and discuss under what conditions one can make 
this  operator upper triangular  by a coordinate change. 
 
Since $J_n\bigl(\lambda(x)\bigr)$ and $L(x)$ are similar,
there exists a matrix-valued  smooth function  $\widehat A(x)$ such that 
\begin{equation} 
\label{eq:L(x)} L(x)= \widehat A^{-1}(x) J_n\bigl(\lambda(x)\bigr) \widehat A (x). 
\end{equation} 
Without loss of generality we may assume that $A(0,\dots,0)= \textrm{\bf 1} = \operatorname{diag}(1,\dots, 1)$.

Consider now the linearization of $L(x)$,  near the point $0=(0,\dots,0)$, 
i.e., take the linear approximation of  $L(x)$ given by \eqref{eq:L(x)}. 
Clearly, up to second order terms, we have  \begin{equation}
 \label{eq:table}   
\begin{array}{ccc}
   \widehat A(x)   &\simeq&  {\textrm{ \bf 1}} +  A \\
     \widehat A^{-1}(x) &\simeq&{\textrm{\bf 1}} - A \\
    {J_n\bigl(\lambda(x)\bigr)} & \simeq & J_n( \lambda(0))+ \Lambda\cdot \textrm{\bf 1} \, , 
 \end{array}\end{equation} 
where $A(x)$ is a certain matrix whose entries are linear functions in local coordinates $x_1,\dots,x_n$: 
$$
A^i_j = \sum_{k} a^i_{j;k} x_k \ , \quad a^i_{j;k}\in\mathbb R,
$$
and $\Lambda= \sum_n\lambda_{k}x_k   $ with some constants $\lambda_{k}$. 

 Substituting \eqref{eq:table} in \eqref{eq:L(x)}, we obtain, up to second order  terms, 
\begin{eqnarray} 
\label{eqarray1} 
L(x) &\simeq & (\textrm{\bf 1} - A(x)  ) (J_n( \lambda(0))+ \Lambda\cdot \textrm{\bf 1})  (\textrm{\bf 1} +  A(x)  )\\ &\simeq&  J_n( \lambda(0))+ \Lambda(x)\cdot \textrm{\bf 1}- A(x) J_n(\lambda(0)) + J_n(\lambda(0))  A(x).   \label{eqarray2} 
\end{eqnarray}
Next, recall that the components of the generalized Nijenhuis torsion  $\mathcal T_L^{(k)}$
of any level $k\ge 1$ (and in particular, the components $\mathcal H^i_{jk}$ of the Haantjes torsion $\mathcal H_L$) are algebraic expressions in the components of $L$ and  their first derivatives; moreover, the derivatives come linearly.  At the point $0=(0,\dots, 0)$, the components of $L$ and their first derivatives  coincide with those of \eqref{eqarray2}. Therefore, the conditions  
${\mathcal H}^i_{jk} {}_{| x =0}=0$  and $T^i_{jk} {}_{| x=0}=0$  are explicit systems of linear equations 
on $a^i_{j;k}$ and $\lambda_k$ viewed as unknowns, whose coefficients  may a priori depend on $\lambda(0)$. It is known though, 
see e.g.  \cite[\S II]{Bo2006}, that
the Haantjes torsions  of $L$ and of $L- \lambda(x) \cdot\textrm{\bf 1}$ coincide. Therefore, $\lambda(0)$ and $\lambda_k$ do not appear in 
the condition 
${\mathcal H}^i_{jk} {}_{| x =0}=0$ related to the proof of Theorem \ref{thm:dim3a}. For the same reason,  the components of $T^i_{jk}$ do not depend on $\lambda(x)$ either as in the definition \eqref{eq:tensorT} we use only the Haantjes torsion and traceless part $\widehat L$ of $L$.   
Thus, the systems of equations  ${\mathcal H}^i_{jk} {}_{| x =0}=0$  and $T^i_{jk} {}_{| x=0}=0$ involve only $a^\alpha_{\beta;\gamma}$ as unknown variables with 
some fixed constant coefficients.

Next, assume that $L(x)$ can be reduced to a triangular form.  It is equivalent to the condition that for any $k$ the distribution  $\operatorname{Ker}(A-\lambda\, \textrm{\bf 1})^k$ is integrable.   The distribution $\operatorname{Ker}(A-\lambda \,\textrm{\bf 1})^k$
is clearly generated by the  vectors $\widehat A^{-1}\partial_{x_1}, \dots, \widehat A^{-1} \partial_{x_k}$, and the integrability condition of this distribution is just the condition that 
\begin{equation}
\label{eq:cond}
[\widehat A^{-1}\partial_{x_i},  \widehat A^{-1} \partial_{x_j}]
\in \textrm{span}_{k\le \max(i,j)}{\left(\widehat A^{-1}\partial_{x_{k}}\right)}.
\end{equation} 
Since the commutator contains only the first derivatives of  vector fields, the condition \eqref{eq:cond}, evaluated at the point $0$, is the following 
linear system of equations on  $a^i_{j;k}$:
\begin{equation} \label{eq:cond3} 
0=[(\textrm{\bf 1}- A)\partial_{x_i}, (\textrm{\bf 1}- A)\partial_{x_j}]_{|x=0}^k = -a^k_{j; i}+ a^k_{i;j}  \  \  \textrm{for $1\le i<j<k\le n$}.
\end{equation}
(The system can be of course immediately solved).

Thus,  we have two systems of equations on $a^k_ {i;j}$ for each dimension, 3 and 4.  If we check that these systems of equations are {\it algebraically}
equivalent, in the sense that any solution of the first is a solution of the second and vice versa, we will show that for $x=0$ the vanishing of $\mathcal H_L$ in dimension 3 and of $T$ in dimension 4 implies  the integrability condition for the distributions 
$\operatorname{Ker}(A-\lambda \,\textrm{\bf 1})^k$ and vice versa. As  there is no essential difference between the point $x=0$ and any other point, and as
the  conditions controlled by the systems  are geometric and do not depend on the choice of coordinate system, 
vanishing of the Haantjes torsion  at every point would imply, in dimension 3, the integrability  
$\operatorname{Ker}(A-\lambda \,\textrm{\bf 1})^k$ and hence upper triangularisibility of $L$, and vice versa. Similarly, in dimension $4$
the condition $T=0$ fulfilled 
at every point would imply upper triangularisibility of $L$, and vice versa.

Now, it is easy to check by direct calculations that the  two systems of equations appeared in the context of Theorem \ref{thm:dim3a} are indeed algebraically equivalent.  In dimension $n=3$, by direct calculations of the Haantjes torsion for the operator at $x=0$, we see that the 
only potentially nonzero components  of $\mathcal H$ are  
$$
\mathcal{H}^1_{2,3}= - \mathcal{H}^1_{3,2} = 3\,a^3_{{1;2}}-3\,a^3_{{2;1}}.
$$ 
Hence, the vanishing of $\mathcal H$ is equivalent to   \eqref{eq:cond3}. This proves  Theorem \ref{thm:dim3a}.

Similarly, in dimension $4$, by direct calculation of $T$ for the operator \eqref{eqarray2}, we see that  the non-zero components  $T^i_{jk}$ are: 
$$
      \begin{array}{rcccr} T^1_{ 2 4}& =&  -2 a^4_{1 ;2} + 2 a^{4}_{ 2
;1}&=& T^1_{4 2}\\
       T^1_{33}&=& 4 a^4_{1; 2} - 4 a^4_{ 2; 1}&&\\
T^1_{3 4}&=&  -a^3_{1; 2} - 2 a^4_{3; 1} + a^3_{2; 1} + 2 a^4_{1;3} && \\
T^1_{43} & =& \ \ a^3_{1; 2} - 2 a^4_{3; 1} - a^3_{2; 1} + 2 a^4_{1;3}\\
      T^1_{44} &=&-4 a^4_{ 3;2}+ 4 a^4_{ 2;3} && \\
    T^2_{34}& =&  -a^4_{1;2} + a^{4}_{ 2; 1} &=& T^2_{43}.
\end{array}
$$

Equating $T$ 
 to zero gives the system algebraically equivalent to \eqref{eq:cond3}.

\begin{remark}{\rm
In the context of the above proof,  the existence of nilpotent operators $L$ as in   
    Example \ref{ex:1}, which are not reducible to a triangular form although   
    their Haantjes torsion $\mathcal T^{(3)}_L$ of level 3 vanishes, can be explained as follows:  the system \eqref{eq:cond3} responsible for 
    upper triangularisibility in dimension 4 gives four linearly independent conditions  on  $a^i_{j;k}$, while the vanishing of  $\mathcal T^{(3)}_L$ gives 
    only two. 
}\end{remark}

\begin{remark} \label{rem:1}{\rm
    In dimension $4$, vanishing of  the Haantjes torsion {$\mathcal T^{(2)}_L= \mathcal{H}^i_{jk}$}  of the  operator \eqref{eqarray2}
    gives 6 
    independent linear equations on $a^i_{j;k}$ at the point $0$, while   
    \eqref{eq:cond3}  contains only 4 independent linear equations.  Indeed, as mentioned above and observed in \cite{TT22},  
   in dimension 4,  not every upper triangular operator has zero Haantjes torsion.  Similarly, not every operator in a 
    strictly upper triangular form has zero Haantjes torsion. On the other hand, vanishing of the Haantjes torsion, viewed as a linear system 
of equations for $a^i_{j;k}$, implies  \eqref{eq:cond3} and therefore local upper triagonalisibility.  One can show that the  later statement holds true in all dimensions and a natural generalization of this statement, which will be published elsewhere,  holds for arbitrary gl-regular operators with real eigenvalues; that is, one can show that  vanishing of the Haantjes torsion is a  sufficient condition for local Jordan--Chevalley decomposition of gl-regular operator fields  with real eigenvalues}.
\end{remark}

\subsection{How did we find the tensor  \texorpdfstring{$T$}{T} from Theorem \ref{thm:2}?}\label{how}

The methods used in the proof of Theorems \ref{thm:dim3a}, \ref{thm:2} allow one to check, by relatively simple 
 calculations, whether a tensor field constructed by $L$ and such that its components are linear in the first 
 derivatives of  $L^i_j$, is `responsible' for the existence of a coordinate system in which $L$ is upper triangular. Let us explain  how we  found the tensor field $T^i_{jk}$. 

Consider the tensor fields of type (1,2) constructed algebraically from $L$ and its Nijenhuis torsion $\mathcal N= \mathcal{T}_L$:\begin{equation}
    \mathcal{N}= \mathcal{N}^i_{jk}, \, L\mathcal{N}= \mathcal{N}^s_{jk} L^i_s,    \,  \mathcal{N} L= \mathcal{N}^i_{sk} L^s_j, \,  (\mathcal{N}L)^\top= \mathcal{N}^i_{js} L^s_k, \,  L\mathcal{N}L= \mathcal{N}^s_{rk} L_s ^iL^r_j, \ \dots\, . \label{eq:4}
 \end{equation}
Then, take linear combination of these tensor fields, with unknown coefficients $c_1,\dots, c_m$. 
Substitution of  $L$ given by \eqref{eq:cond3} and $x=0$ gives a system of linear equations on $a^i_{j;k}$, $\lambda_k$, whose coefficients depend on $\lambda(0)$.  The system  of  equations so obtained can be solved with respect to $c_1,\dots,c_m$. For any solution, the corresponding tensor vanishes if $L$ can be put in the upper diagonal form; in other words, any choice of a 
solution of the system gives us a necessary condition for the existence of a coordinate system such that $L$ is upper triangular.
Finally, one needs to  check whether  for a generic choice of a solution, the vanishing of the corresponding tensor at $x=0$ is equivalent to \eqref{eq:cond3}. This is expected, if we pick  sufficiently many tensors of form \eqref{eq:4}.

 \section{Proof of the  Tempesta--Tondo Conjecture}
 
We consider two operators  $L$, $M$ given by strictly upper triangular matrices in a local coordinate system $x_1,\dots,x_n$.
From the definition  \eqref{haam} of the generalized Haantjes torsion of level $(n -1)$, we see that  
$\mathcal{H}^{(n-1)}_{L,M}(\xi ,\eta)$  is  a linear combination of terms of the form 
\begin{equation}\label{eq:power}
   L^{\alpha_1}M^{\beta_1} [ L^{\alpha_2} M^{\beta_2} \xi ,  L^{\alpha_3} M^{\beta_3} \eta ] \  \ \textrm{with $\sum_{i=1}^3 (\alpha_i+ \beta_i)
   = 2n -2.  $ }
\end{equation}
Here $\alpha_i, \beta_i\in \mathbb{N}\cup \{0\}$, \  $L^\alpha \zeta$ is the product of $\alpha$ copies of the matrix $L$ and $\zeta$ treated as a column-vector, and the square brackets denote the Lie bracket of vector fields. 

As the Haantjes torsion of level $(n -1)$ is a tensorial object, it is sufficient to check  vanishing of \eqref{eq:power} for the basis vector fileds 
$\partial_{x_i}$.  Since $L$ and $M$  are strictly upper triangular, then 
$L\partial_{x_i} \in \textrm{span}_{j<i}(\partial_{x_j})$ and $M\partial_{x_i} \in \textrm{span}_{j<i}(\partial_{x_j})$. In particular, $L^{n}\partial_{x_i}=0$. 

Moreover, since the distribution  $\textrm{span}_{k\le i}(\partial_{x_k})$
is integrable for any $i$, we have that for any vectors $\xi   \in  \textrm{span}_{k\le i}(\partial_{x_k})$ 
and  $\eta  \in  \textrm{span}_{k\le j}(\partial_{x_k})$      
$$
[ \xi  , \eta]\in   \textrm{span}_{k \le \operatorname{max}(i,j)}(\partial_{x_k}), \quad i,j=1,\dots,n.
$$
This implies that  the terms of the form \eqref{eq:power}
such that  $\alpha_1 + \beta_1\ge n$,  $\alpha_2 + \beta_2\ge n$ or  $\alpha_3 + \beta_3\ge n$  automatically vanish.

Assume  $\alpha_2 +\beta_2 \ge \alpha_3 + \beta_3$.  Let us show that the term  
$L^{\alpha_1}M^{\beta_1}[ L^{\alpha_2} M^{\beta_2} \xi ,  L^{\alpha_3} M^{\beta_3} \eta ] $ vanishes  unless  $\alpha_2+ \beta_2= n-1$.
Indeed, if  $\alpha_2 +\beta_2\ge n$, then 
$L^{\alpha_2} M^{\beta_2} \xi =0$ and the statement follows.  Next we observe that  
$$ 
L^{\alpha_1}M^{\beta_1}[ L^{\alpha_2} M^{\beta_2} \xi ,  L^{\alpha_3} M^{\beta_3} \eta ] \in \textrm{span}_{k \le n-\alpha_3 - \beta_3-\alpha_1 - \beta_1 }(\partial_{x_k}).
$$
Thus, if the term does not vanish, then $1\le n-\alpha_3 - \beta_3-\alpha_1 - \beta_1$. Combining this with  \eqref{eq:power}, we obtain $\alpha_2 +\beta_2\ge n-1$. 

Similarly, if $\alpha_3 +\beta_3 \ge \alpha_2 + \beta_2$, the term  
$L^{\alpha_1}M^{\beta_1}[ L^{\alpha_2} M^{\beta_2} \xi ,  L^{\alpha_3} M^{\beta_3} \eta ]$ vanishes unless $\alpha_3 + \beta_3 =n-1$.

Thus, we may assume that either $\alpha_2 + \beta_2= n-1$ or $\alpha_3 + \beta_3= n-1$. Without loss of generality, we consider the first case, 
$\alpha_2 + \beta_2= n-1$. Then in view of \eqref{eq:power},  we have 
\begin{equation}
\label{eq:power2}
\alpha_1+\beta_1+\alpha_3 + \beta_3 = n-1.
\end{equation}

Now, observe that  if $i\le n-1$ or $j \le n-1$, then 
$$
L^{\alpha_1}M^{\beta_1}[ L^{\alpha_2} M^{\beta_2} \partial_{x_i} ,  L^{\alpha_3} M^{\beta_3} \partial_{x_j} ] =0.
$$
Indeed, since $\alpha_2 + \beta_2= n-1$, then $L^{\alpha_2} M^{\beta_2} \partial_{x_i}$ vanishes unless $i=n$.  If $i=n$, then $L^{\alpha_2} M^{\beta_2} \partial_{x_i}$ is proportional to $\partial_{x_1}$ and therefore
\begin{equation} 
\label{eq:rem} 
L^{\alpha_1}M^{\beta_1}[ L^{\alpha_2} M^{\beta_2} \partial_{x_i} ,  L^{\alpha_3} M^{\beta_3} \partial_{x_j} ]\in 
\textrm{span}_{k \le j-\alpha_3 - \beta_3-\alpha_1 - \beta_1 }(\partial_{x_k}). 
\end{equation}   
In view of \eqref{eq:power2},  we see that 
$L^{\alpha_1}M^{\beta_1}[ L^{\alpha_2} M^{\beta_2} \partial_{x_i} ,  L^{\alpha_3} M^{\beta_3} \partial_{x_j} ]=0$ unless $j= n$.

Thus, we have proved that  $\mathcal{H}^{(n-1)}_{L,M}(\partial_{x_i},\partial_{x_j})= 0$, if $i<n$ or $j<n$. To complete the proof of the Tempesta--Tondo conjecture, it remains to check that 
$\mathcal{H}^{(n-1)}_{L,M}(\partial_{x_n},\partial_{x_n})= 0$. But this is trivial as $\mathcal{H}^{(n-1)}_{L,M}$ is  skew-symmetric.

\printbibliography

\end{document}